\documentclass[12pt,reqno]{amsart}
\usepackage{amsmath, amsthm, amscd, amsfonts, amssymb, graphicx}%, color}

\setbox0=\hbox{$+$}
\newdimen\plusheight
\plusheight=\ht0
\def\+{\;\lower\plusheight\hbox{$+$}\;}

\setbox0=\hbox{$-$}
\newdimen\minusheight
\minusheight=\ht0
\def\-{\;\lower\minusheight\hbox{$-$}\;}

\setbox0=\hbox{$\cdots$}
\newdimen\cdotsheight
\cdotsheight=\plusheight%\ht0
\def\cds{\lower\cdotsheight\hbox{$\cdots$}}

\renewcommand{\(}{\left\(}
\renewcommand{\)}{\right\)}

\renewcommand{\pmod}[1]{\,(\textup{mod}\,#1)}
\numberwithin{equation}{section}
\theoremstyle{plain}
\newtheorem{theorem}{Theorem}[section]
\newtheorem{lemma}[theorem]{Lemma}

\newtheorem{remark}[theorem]{Remark}

\newtheorem*{theorem*}{Theorem}

\begin{document}

\setcounter{page}{1}

 \centerline{\bf \large Integer Partitions With Restricted Distinct Parts} \vskip3mm

 \begin{center}
 {\bf Rinchin Drema and Nipen Saikia$^\ast$}\\
 {Department of Mathematics, Rajiv Gandhi University,\\ Rono Hills, Doimukh-791112, 
 Arunachal Pradesh, India.
 } { e-mail: {dremarinchin3@gmail.com; nipennak@yahoo.com}}\\
 {\bf $^\ast$ Corresponding author}.
  
 \end{center} 
 
 \vskip 1cm
 \centerline{\bf Abstract}
  For any positive integers $s$ and $t$, let $Q_{t}^{s}(n)$ denotes the number of partitions of a positive integer $n$ into distinct parts such that no part is congruent to $s$ or $t-s$ modulo $t$. We prove some Ramanujan-type congruences for $Q_{t}^{s}(n)$ for some particular values of $s$ and $t$  by employing $q$-series and theta function identities.
 \vskip 3mm
   
  \vskip 0.5cm
  

\section{Introduction}
The partition of a non-negative integer $n$ is a non-increasing sequence of positive integers  called parts, where the sum of all the parts is equal to $n$.
 For example,  the positive integer $n=5$ has seven partitions, viz,$$ 5, \quad 4+1, \quad 3+2, \quad 3+1+1, \quad 2+2+1, \quad2+1+1+1,\quad 1+1+1+1+1.$$
The number of partition of a non-negative integer $n$ is generally denoted by $p(n)$ with $p(0)=1$. Thus,  $p(5)=7$.

The generating function for the partition function $p(n)$ due to Euler \cite{eu} is given by 
\begin{equation}\label{pngeni}\sum_{n\ge 0}p(n)q^n=\dfrac{1}{(q;q)_{\infty}},\end{equation}where for any complex number $a$ and $|q|<1$ , 
$$(a;q)_{\infty}=\prod_{n\ge 0}(1-aq^n).$$
Throughout this paper, we will use the notation, for any positive integer $k$, 
\begin{equation}\label{fk}f_k:=(q^k;q^k)_{\infty}.\end{equation} Euler also considered two restricted partition functions, the number of partitions of $n$ into distinct parts denoted by $p_d ( n)$  and the number of partitions of $n$ into odd parts denoted by $ p_o (n)$ and proved that famous theorem, $p_d (n)=p_o (n)$. In terms of generating functions, the theorem is stated as \begin{equation}\sum_{n\ge 0}p_{o}(n)q^n=\frac{1}{(q;q^2)_\infty}=(-q;q)_\infty=\sum_{n\ge 0}p_{d}(n)q^n.\end{equation}

Ramanujan \cite{sr} established following congruences for the partition function $p(n)$:
$$p(5n+4) \equiv0\pmod5,$$
$$p(7n+5)\equiv0 \pmod7,$$ and
$$p(11n+6)\equiv0\pmod {11}.$$ 
A large numbers of mathematicians are inspired by Ramanujan's congruences for $p(n)$. This led mathematicians to work in search of congruence properties for other partition functions and in this process, many restricted partition functions are also introduced and studies are  presently going on for their arithmetic properties. For instance, see \cite{SS, kw}.
In this paper, we consider the new restricted partition function denoted by  $Q_{t}^{s}(n)$ which counts the number of distinct partitions of $n$ such that no part is congruent to $s$ or $t-s$ modulo $t$. The generating function of  $Q_{t}^{s}(n)$ is given by
 \begin{equation} \label{r1}
\sum_{n\ge 0}Q_{t}^{s}(n)q^n=\dfrac{(-q;q)_{\infty}}{(-q^s;q^t)_{\infty}(-q^{t-s};q^t)_\infty}.\end{equation}\\
Applying elementary $q$-operations, it is easily seen that
$$
\sum_{n\ge 0}Q_{t}^{s}(n)q^n =\dfrac{(-q;q)_{\infty}(q;q)_{\infty}(q^t;q^t)_{\infty}}{(q;q)_{\infty}(-q^s;q^t)_{\infty}(-q^{t-s};q^t)_{\infty}(q^t;q^t)_{\infty}}$$
\begin{equation} \label{z2}\hskip 1.5cm=\dfrac{(q^2;q^2)_{\infty}(q^t;q^t)_{\infty}}{(q;q)_{\infty}f(q^s,q^{t-s})} =\dfrac{f_2 f_t}{f_1 f(q^s,q^{t-s})},\end{equation}
where $f(a,b)$ \cite[Entry 19, p.35]{bc3} given by
 \begin{equation}\label{h2}
 f(a,b)=(-a;ab)_{\infty}(-b;ab)_{\infty}(ab;ab)_{\infty}=\sum_{n=-\infty}^\infty a^{n(n+1)/2}b^{n(n-1)/2},~  |ab|<1\end{equation}
 is the Ramanujan's general theta-function.

The purpose of this paper is to establish some Ramanujan-type congruences for $Q_{t}^{s}(n)$ in the spirit of Ramanujan and its connections with some other partition fucntions. To prove our results, we will employ some $q$-series and theta-function identities which is also in the spirit of Ramanujan. We also give suitable examples to verify our congruences.

\section{Preliminaries}

This section is devoted to record some theta-function identities for ready refernces in this paper. 

Three useful special cases of $f(a, b)$ in \eqref{h2} are  the theta-functions $\phi(q)$, $\psi(q)$ and $f(-q)$  \cite[p. 36, Entry 22(iii)]{bc3} given by  
 \begin{equation}\label{e0}
\phi(q):=f(q, q)=\sum_{n\ge 0}q^{n^2}=(-q;q^2)^{2}_{\infty}(q^2:q^2)_{\infty}=\dfrac{{f
_2}^5}{{f_1}^2{f_4}^2}, 
 \end{equation} 
\begin{equation}\label{e1}
\psi(q):=f(q,q^3)=\sum_{n\geq 0}q^{n(n+1)/{2}}=\dfrac{(q^2;q^2)_{\infty}}{(q:q^2)_{\infty}}=\dfrac{{f
_2}^2}{f_1},
\end{equation} and
\begin{equation}\label{e2} 
 f(-q):=f(-q, -q^2)=\sum_{n\geq 0}(-1)^n q^{n(3n-1)/{2}}=(q;q)_{\infty}=f_1 .\end{equation}
 Ramanujan also defined the function $\chi(q)$ as 
 \begin{equation}\chi(q)=(-q;q^2)_\infty. 
 \end{equation}
Employing elementary $q$-operations, it is easily seen that \begin{equation}\label{e3}
\phi(-q)=\dfrac{(q,q)^{2}_{\infty}}{(q^2;q^2)_{\infty}}=\dfrac{f_{1}^2}{f_2},
\end{equation}
\begin{equation}
\psi(-q)=\dfrac{(q;q)_{\infty}(q^4;q^4)_{\infty}}{(q^2;q^2)_{\infty}}=\dfrac{f_1 f_4}{f_2},
\end{equation} and
 \begin{equation} 
\chi(q)=\dfrac{(q^2;q^2)_{\infty}^{2}}{(q;q)_{\infty}(q^4;q^4)_{\infty}}=\dfrac{{f_2}^2}{f_1 f_4}.
\end{equation}
 \begin{lemma}\label{lem2}\cite[p. 350, Eq.(2.3)]{bc3} We have
 \begin{equation}\label{t8}
 f(q,q^2)=\dfrac{\phi(-q^3)}{\chi(-q)},\end{equation}
 where,
 \begin{equation}\label{lo}
 \chi(-q)=\dfrac{(q;q)_{\infty}}{(q^2;q^2)_{\infty}}=\dfrac{f_1}{f_2}.
 \end{equation}
 \end{lemma}
 
   \begin{lemma}\cite [p. 51, Example(v)]{bc3} We have
   \begin{equation}\label{t9}
  f(q,q^5)={\psi(-q^3)}{\chi(q)}.\end{equation}
  \end{lemma}
  
  \begin{lemma} \cite[Theorem 2.2]{cg} For any prime $p\geq 5$, we have
   $$
 f_1=\sum_{\substack{k={-(p-1)/2} \\ k \ne {(\pm p-1)/6}}}^{(p-1)/2}(-1)^{k}{\it q}^{(3m^2+m)/2} f\left({-\it q}^{{(3p^2+(6k+1)p)/2}},{-\it q}^{(3p^2-(6k+1)p)/2}\right)$$\begin{equation}\label{u10} +(-1)^{{(\pm p-1)/6}}{\it q}^{(p^2-1)/24}f_{p^2}.\end{equation}where
 \begin{equation*}
 \dfrac{\pm p-1}{6}
 = \left\{
         \begin{array}{ll}
          \dfrac{(p-1)}{6},   
            & if~ p \equiv 1\pmod 6 \\
            \dfrac{(-p-1)}{6} ,
            & if~ p \equiv -1\pmod 6. 
         \end{array}
     \right.
 \end{equation*}Furthermore,  if $$\dfrac{-(p-1)}{2}\leq k \leq\dfrac{(p-1)}{2}~and ~k \neq \dfrac{(\pm p-1)}{6},$$~ then\\
 $$\dfrac{3k^2+k}{2}\not\equiv \dfrac{p^2-1}{24} \pmod p.$$
 \end{lemma}
 \begin{lemma}\cite[Theorem 2.1]{cg} For any odd prime $p$, we have
 \begin{equation}\label{p8}
 \psi(q)= \sum_{m=0}^{(p-3)/2}{\it q}^{(m^2+m)/2} f\left( {\it q}^{(p^2+(2m+1)p)/2},{\it q}^{(p^2-(2m+1)p)/2}\right)+{\it q}^{(p^2-1)/8}\psi(q^{p^2}).\end{equation} 
 Furthermore, $\dfrac{m^2+m}{2}\not \equiv \dfrac{p^2-1}{8} \pmod p$ for $0\leq m\leq\dfrac{p-3}{2}.$
 
 \end{lemma}
 \begin{lemma}\cite[p. 303, Entry 17(v)]{bc3} We have
 \begin{equation}\label{u7}
 f_1=f_{49}\left(\dfrac{B(q^7)}{C(q^7)}-q \dfrac{A(q^7)}{B(q^7)}-q^2+q^5\dfrac{C(q^7)}{A(q^7)}\right),
 \end{equation}
 where $A(q)=f(-q^3,-q^4),  B(q)=f(-q^2,-q^5)~and ~C(q)=f(-q,-q^6).$
 \end{lemma}
 
 \begin{lemma}\cite{sr} We have \begin{equation} \label{g1}
 f_1=f_{25}\left(R(q^5)-q-q^2R(q^5)^{-1}\right),
 \end{equation} where $$R(q)=\dfrac{(q^2;q^5)_{\infty}(q^3;q^5)_{\infty}}{(q;q^5)_{\infty}(q^4;q^5)_{\infty}}.$$
 \end{lemma}
 \begin{lemma}\label{g2}\cite[ eq.(21.3.3)]{HS4} We have
 \begin{equation}\label{y9}
  f_1^3=a(q^3)-3qf_9^3,
 \end{equation} 
 where $a(q)=\dfrac{f_2^6 f_3}{f_1^2 f_6^2} +3q\dfrac{f_1^2 f_6^6}{f_2^2 f_3^3}.$ \end{lemma}
 
 We end this section by recording some congruence properties.
Employing binomial theorem in \eqref{fk}, one can easily see that, 
  for all positive integer $r$ and $m$,
  \begin{equation}\label{t7}
 { f_r}^{2m}\equiv {f_{2r}}^m\pmod2\end{equation} and
 \begin{equation}\label{v7}
  { f_r}^{4m}\equiv {f_{2r}}^{2m}\pmod4.\end{equation} 
 
 \begin{lemma} We have
 \begin{equation}\label{jk9}
 f(q,q^2)\equiv f_1 \pmod 2. 
\end{equation}
\begin{proof}
Using \eqref{e3} in \eqref{t8} and employing \eqref{t7}, we arrive at the desired result.
\end{proof}
  \end{lemma}

\section{Congruences of $Q_{t}^{s}(n)$}
\begin{theorem}
For any prime $p \geq 5$ with $\Big(\dfrac{-2\alpha}{p}\Big)=-1$, $ 1\leq j\leq p-1$, $\alpha\in N $ and $ \beta \geq 0,$  we have
 \begin{equation}\label{n20}
 \sum_{n\geq0} Q_{2\alpha}^{\alpha}\left(p^{2\beta}n+\dfrac{(2\alpha+1)(p^{2\beta} -1)}{24} \right) \equiv f_{1} f_{2\alpha} \pmod 2,
 \end{equation} 
 and 
  \begin{equation}\label{n21}
 \sum_{n\geq0} Q_{2\alpha}^{\alpha}\left(p^{2\beta +1}(pn+j)+\dfrac{(2\alpha+1)(p^{2\beta+2} -1)}{24} \right) \equiv 0 \pmod 2.
 \end{equation}
\end{theorem}
\begin{proof}Our proof relies on mathematical induction.
Setting $t=2\alpha$ and $s=\alpha$ in \eqref{z2}, we obtain
\begin{equation}\label{i2}
\sum_{n\geq0} Q_{2\alpha}^{\alpha}(n)q^n=\dfrac{f_2 f_{2\alpha}}{f_1 f(q^{\alpha},q^{\alpha})}.
\end{equation}
Now using \eqref{e0} in \eqref{i2}, we obtain
\begin{equation}\label{i3}
\sum_{n\geq0} Q_{2\alpha}^{\alpha}(n)q^n =\dfrac{f_2 f_{2\alpha}f_{\alpha}^2{f_{4\alpha}^2}}{f_1 f_{2\alpha}^5}.\end{equation}
Using \eqref{t7} in \eqref{i3}, we obtain 
\begin{equation}\label{i4}
\sum_{n\geq0} Q_{2\alpha}^{\alpha}(n)q^n \equiv f_1 f_{2\alpha} \pmod 2,\end{equation} 
which is the  case $\beta=0$ of \eqref{n20}. Assume that \eqref{n20} holds for some $\beta \geq 0.$ 

Employing  \eqref{u10} in \eqref{n20}, we obtain   
\begin{multline}\label{i13}
 \sum_{n\geq0} Q_{2\alpha}^{\alpha}\left( p^{2\beta}n+\dfrac{(2\alpha+1)(p^{2\beta} -1)}{24} \right) \\\equiv  \Big[\sum_{\substack{k={-(p-1)/2} \\ k \ne {(\pm p-1)/6}}}^{(p-1)/2}(-1)^{k}{\it q}^{(3k^2+k)/2} f\left( {-\it q}^{{(3p^2+(6k+1)p)/2}},{-\it q}^{(3p^2-(6k+1)p)/2}\right)   
\\ +(-1)^{{(\pm p-1)/6}}{\it q}^{(p^2-1)/24}f_{p^2} \Big]\\\times \Big[ \sum_{\substack{ m={-(p-1)/2} \\ m \ne {(\pm p-1)/6}}}^{(p-1)/2}(-1)^{m}{\it q}^{\alpha(3m^2+m)} f\left( {-\it q}^{{2\alpha(3p^2+(6m+1)p)/2}},{-\it q}^{2\alpha(3p^2-(6m+1)p)/2}\right)\\   
 +(-1)^{{(\pm p-1)/6}}{\it q}^{2\alpha(p^2-1)/24}f_{2\alpha p^2} \Big] \pmod 2. \end{multline}
Consider the congruence
$$\dfrac{3k^2+k}{2}+{\alpha(3m^2+m)} \equiv \dfrac{(2\alpha+1)(p^2-1)}{24} \pmod p,$$
which is equivalent to 
$$ (6k+1)^2+2\alpha(6m+1)^2 \equiv 0 \pmod p.$$
For $\Big(\dfrac{-2\alpha}{p}\Big)=-1,$ the  above congruence  has unique solution $k=m=\dfrac{(\pm p-1)}{6}.$ Therefore, extracting the term involving $q^{pn+{{(2\alpha+1)(p^2-1)/24}}}$ from both sides of
\eqref{i13}, dividing  by $q^{{(2\alpha+1)(p^2-1)/24}}$ and then replacing $q^{2p}$ by $q,$ we obtain
\begin{equation}\label{i14}
 \sum_{n\geq0} Q_{2\alpha}^{\alpha}\left(p^{2\beta+1}n +\dfrac{(2\alpha+1)(p^{2\beta+2} -1)}{24} \right) \equiv f_{p} f_{2\alpha p} \pmod 2.
 \end{equation}  
Extracting the terms involving $q^{pn}$ from \eqref{i14} and replacing $q^p $ by $q$, we obtain
\begin{equation}\label{d10}
 \sum_{n\geq0} Q_{2\alpha}^{\alpha}\left(p^{2\beta+2}n+\dfrac{(2\alpha+1)(p^{2\beta+2} -1)}{24} \right) \equiv f_{1} f_{2\alpha} \pmod 2,
 \end{equation}
 which is the case $\beta +1$ of \eqref{n20}. So the proof of \eqref{n20} is complete.
 
 Now, extracting the coefficients of $q^{pn+j},$ for $1\leq j \leq p-1,$ from both sides of \eqref{i14}, we arrive at \eqref{n21}.
 \end{proof}
\begin{theorem}
For any prime $p \geq 5$ with $\Big(\dfrac{-\alpha}{p}\Big)=-1$, $ 1\leq j\leq p-1$, $\alpha\in N $ and $ \beta \geq 0,$  we have
 \begin{equation}\label{t10}
 \sum_{n\geq0} Q_{4\alpha}^{\alpha}\left(p^{2\beta}n+\dfrac{(\alpha+1)(p^{2\beta} -1)}{24} \right) \equiv f_{1} f_{\alpha} \pmod 2,
 \end{equation} 
 and  
  \begin{equation}\label{t11}
 \sum_{n\geq0} Q_{4\alpha}^{\alpha}\left(p^{2\beta +1}(pn+j)+\dfrac{(\alpha+1)(p^{2\beta+2} -1)}{24} \right) \equiv 0 \pmod 2.
 \end{equation}
\end{theorem}
\begin{proof} 
Setting $t=4\alpha$ and $s=\alpha$ in \eqref{z2}, we obtain
\begin{equation}\label{r2}
\sum_{n\geq0} Q_{4\alpha}^{\alpha}(n)q^n=\dfrac{f_2 f_{4\alpha}}{f_1 f(q^{\alpha},q^{3\alpha})}.
\end{equation}
Using \eqref{e1} in \eqref{r2}, we obtain
 \begin{equation}\label{r3}
\sum_{n\geq0} Q_{4\alpha}^{\alpha}(n)q^n=\dfrac{f_2 f_{4\alpha}}{f_1 \psi(q^{\alpha})}.
\end{equation}
Employing \eqref{e1} in \eqref{r3} and then using  \eqref{t7}, we obtain
\begin{equation}\label{r4}
\sum_{n\geq0} Q_{4\alpha}^{\alpha}(n)q^n \equiv f_1 f_{\alpha} \pmod 2,\end{equation} 
which is the case $\beta=0$ of \eqref{t10}. Assume that \eqref{t10} holds for some $\beta \geq 0$. Employing \eqref{u10} in \eqref{t10}, we obtain  
\begin{multline}\label{q13}
 \sum_{n\geq 0} Q_{4\alpha}^{\alpha}\left( p^{2\beta}n+\dfrac{(\alpha+1)(p^{2\beta} -1)}{24} \right)\\ \equiv \Big[\sum_{\substack{ k={-(p-1)/2} \\ k \ne {(\pm p-1)/6}}}^{(p-1)/2}(-1)^{k}{\it q}^{(3k^2+k)/2} f\left( {-\it q}^{{(3p^2+(6k+1)p)/2}},{-\it q}^{(3p^2-(6k+1)p)/2}\right)\\   
 +(-1)^{{(\pm p-1)/6}}{\it q}^{(p^2-1)/24}f_{p^2}\Big]\\ \times\Big[ \sum_{\substack{ m={-(p-1)/2} \\ m \ne {(\pm p-1)/6}}}^{(p-1)/2}(-1)^{m}{\it q}^{\alpha(3m^2+m)/2} f\left( {-\it q}^{{\alpha(3p^2+(6m+1)p)/2}},{-\it q}^{\alpha(3p^2-(6m+1)p)/2}\right)\\   
 +(-1)^{{(\pm p-1)/6}}{\it q}^{\alpha(p^2-1)/24}f_{\alpha p^2}\Big]  \pmod 2. \end{multline}
Consider the congruence 
$$\dfrac{3k^2+k}{2}+\dfrac{\alpha(3m^2+m)}{2} \equiv \dfrac{(\alpha+1)(p^2-1)}{24} \pmod p,$$
which is equivalent to 
$$ (6k+1)^2+\alpha(6m+1)^2 \equiv 0 \pmod p.$$
For $\Big(\dfrac{-\alpha}{p}\Big)=-1,$ the  above congruence  has unique solution $k=m=\dfrac{(\pm p-1)}{6}.$ So extracting the term involving $q^{pn+{(\alpha+1)(p^2-1)/24}}$ from both sides of
\eqref{q13}, dividing by $q^{{(\alpha+1)(p^2-1)/24}}$ and then replacing $q^p$ by $q$, we obtain
\begin{equation}\label{q14}
 \sum_{n\geq0} Q_{4\alpha}^{\alpha}\left(p^{2\beta+1}n+\dfrac{(\alpha+1)(p^{2\beta+2} -1)}{24} \right) \equiv f_{p} f_{\alpha p} \pmod 2.
 \end{equation}  
 Extracting the terms involving $q^{pn}$ from \eqref{q14} and replacing $q^p$ by $q$, we obtain
\begin{equation}\label{y10}
 \sum_{n\geq 0} Q_{4\alpha}^{\alpha}\left(p^{2\beta+2}n+\dfrac{(\alpha+1)(p^{2\beta+2} -1)}{24} \right) \equiv f_{1} f_{\alpha} \pmod 2,
 \end{equation}
 which is the case $\beta +1$ of  \eqref{t10}. So by principle of mathematical induction, proof of \eqref{t10} is complete.
 
  Now extracting the coefficients of $q^{pn+j},$ for $1\leq j \leq p-1,$ from both sides of \eqref{q14}, we arrive at \eqref{t11}.
 \end{proof}

 \begin{theorem}\label{thm3} For any integer $\alpha\ge1,$
 we have
 \begin{equation}\label{p9}
 \sum_{n\geq0} Q_{5\alpha}^{5}\left(5n+i \right) \equiv 0 \pmod 2,\quad i=3,4.
 \end{equation}
 \begin{proof}
 Setting $t=5\alpha $ and $s=5$ in \eqref{z2}, we obtain
  \begin{equation}\label{y11}
 \sum_{n\geq0} Q_{5\alpha}^{5}(n)q^n =\dfrac{f_2 f_{5\alpha}}{f_1 f(q^5,q^{5\alpha-5})}.\end{equation}
 Using \eqref{t7} in \eqref{y11}, we obtain
  \begin{equation}\label{y12}
 \sum_{n\geq0} Q_{5\alpha}^{5}(n)q^n \equiv \dfrac{f_1 f_{5\alpha}}{ f(q^5,q^{5\alpha-5})} \pmod 2. \end{equation}
 Employing \eqref{g1} in \eqref{y12}, we obtain
    \begin{equation}\label{y13}
 \sum_{n\geq0} Q_{5\alpha}^{5}(n)q^n \equiv \dfrac{ f_{5\alpha}}{ f(q^5,q^{5\alpha-5})}\Big(f_{25}(R(q^5)-q-q^2R(q^5)^{-1})\Big) \pmod 2. \end{equation}
 Extracting the term $q^{5n+i},$ for $i=3, 4$  from \eqref{y13}, we arrive at the desired result.
  \end{proof}
 \end{theorem}
 \begin{remark}
Setting $n=1$ and $\alpha=2$ in Theorem \ref{thm3}, we observe that
$Q_{10}^{5}(8)=4$ with the relevant partitions $8 \quad 7+1 \quad 6+2$ and $ 4+3+1$, and $Q_{10}^5(9)=6$ with relevant partitions $9,\quad 8+1,\quad7+2,\quad 6+3,\quad 6+2+1,$ and $ 4+3+2.$
 \end{remark}
\begin{theorem}\label{thm4} For any integer $\alpha\ge 1,$
  we have
 \begin{equation}\label{f9}
 \sum_{n\geq0} Q_{7\alpha}^{7}\left(7n+i \right) \equiv 0 \pmod 2\quad i=,3,4,6.
 \end{equation}
\end{theorem}
\begin{proof}  
 Setting $t=7\alpha $ and $s=7$ in \eqref{z2}, we have
  \begin{equation}\label{y17}
 \sum_{n\geq0} Q_{7\alpha}^{7}(n)q^n =\dfrac{f_2 f_{7\alpha}}{f_1 f(q^7,q^{7\alpha-7})}.\end{equation}
 Using \eqref{t7} in \eqref{y11}, we obtain
  \begin{equation}\label{y18}
 \sum_{n\geq0} Q_{7\alpha}^{7}(n)q^n \equiv \dfrac{f_1 f_{7\alpha}}{ f(q^7,q^{7\alpha-7})} \pmod 2. \end{equation}
 Employing \eqref{u7} in \eqref{y18}, we obtain
    \begin{equation}\label{g13}
 \sum_{n\geq0} Q_{7\alpha}^{7}(n)q^n \equiv \dfrac{ f_{7\alpha}f_{49}}{ f(q^7,q^{7\alpha-7})}\left(\dfrac{B(q^7)}{C(q^7)}-q \dfrac{A(q^7)}{B(q^7)}-q^2+q^5\dfrac{C(q^7)}{A(q^7)}\right) \pmod 2.\end{equation}
 Since the right hand side of \eqref{g13} contains no terms involving $q^{7n+i}$ for $i=3,4,6$,  extracting the term involving $q^{7n+i}$ for $i=3, 4, 6$  from \eqref{g13}, we arrive at the desired result.  \end{proof}
 
 \begin{remark}
Setting $n=1,$ $\alpha=2$ and $i=3$ in Theorem \ref{thm4} we observe that
$Q_{14}^{7}(10)=8$ with the relevant partitions $10, \quad 9+1, \quad 8+2,\quad 6+4,\quad6+3+1,\quad5+4+1,\quad5+3+2,$ and $4+3+2+1.$ Again, setting $n=1,$ $\alpha=2$ and $i=4$ in Theorem \ref{thm4} we observe that $Q_{14}^{7}(11)=10$ with relevant partitions $11,\quad 10+1,\quad 9+2, \quad 8+3, \quad 8+2+1,\quad 6+5,\quad6+4+1,\quad 6+3+2,\quad 5+4+2,$ and $ 5+3+2+1. $
 \end{remark}

 \begin{theorem} We have 
\begin{equation}\label{l9} 
 Q_{3}^{2}(3n+1)=Q_{3}^{2}(3n+2)= 0,\end{equation} and
\begin{equation}\label{43}
Q_{3}^{2}(2n)=p_{o}(n), \end{equation} where $p_{o}(n)$ counts the number of partitions of $n$ with odd parts only.
\end{theorem}
\begin{proof}Setting $t=3$ and $s=2$ in \eqref{z2}, we obtain
\begin{equation}\label{20}
\sum_{n\geq0} Q_{3}^{2}(n)q^n=\dfrac{f_2 f_{3}}{f_1 f(q^2,q)}=\dfrac{f_2 f_{3}}{f_1 f(q,q^2)}. \end{equation}
 Using Lemma \ref{lem2} in \eqref{20}, we obtain
\begin{equation}\label{21}
\sum_{n\geq0} Q_{3}^{2}(n)q^n=\dfrac{ f_{6}}{f_3}. \end{equation}
Since the right hand side contains no terms involving $q^{3n+1}$ and $q^{3n+2}$, extracting the terms involving  $q^{3n+i}$ for $i=1$ and 2, we easily arrive at \eqref{l9}.

Again, extracting the term $q^{3n}$ from \eqref{21} and replacing $q^3$ by $q$, we obtain
\begin{equation}\label{23}
\sum_{n\geq0} Q_{3}^{2}(3n)q^n=\dfrac{ f_{2}}{f_1}=\dfrac{(q^2;q^2)_{\infty}}{(q;q)_{\infty}}. \end{equation}
Simplifying \eqref{23} by using elementary $q$-operations, we obtain
\begin{equation}\label{24}
\sum_{n\geq0} Q_{3}^{2}(3n)q^n=\dfrac{1}{(q;q^2)_{\infty}}=\sum_{n\geq0} p_o(n)q^n. \end{equation} Equating the coefficients of $q^n$ on both sides of \eqref{24}, we arrive at \eqref{43}.
\end{proof}
\begin{remark}
Setting $n=1,2$ and 3 respectively in result \eqref{l9}, we get $Q_3^2(4)=0, Q_3^2(7)=0~and~ Q_3^2(10)=0, respectively.$ One can also verify that  $ Q_3^2(6)=1$ with relevant partition 6 and $p_o(2)=1$ with relevant partition $1+1.$
\end{remark}
\begin{theorem}
Let $p\geq 5$ be any prime with $\left(\dfrac{-6}{p}\right)=-1$ and $1\leq j \leq p-1.$ Then for any non-negative integers $\beta,$
\begin{equation}\label{g7}
\sum_{n\geq0} Q_4^2 \left(p^{2\beta} n +\dfrac{5}{12}(p^{2\beta}-1)\right)q^n \equiv \psi(q)f_2 \pmod 4, \end{equation} 
and 
\begin{equation}\label{s6}
\sum_{n\geq0} Q_4^2 \left(p^{2\beta+1} (pn+j) +\dfrac{5}{12}(p^{2\beta+2}-1)\right)q^n \equiv 0 \pmod 4. \end{equation} 
\end{theorem}
\begin{proof}
Setting $t=4$ and $s=2$ in \eqref{z2}, we obtain
\begin{equation}\label{x6}
\sum_{n\geq0} Q_4^2(n)q^n=\dfrac{f_2 f_4}{f_1 f(q^2,q^2)}.
\end{equation}
Using \eqref{e0} in \eqref{x6}, we obtain
\begin{equation}\label{b8}
\sum_{n\geq0} Q_4^2(n)q^n=\dfrac{{f_2}^3 {f_8}^2}{f_1 {f_4}^4}.
\end{equation}
Employing \eqref{v7} in \eqref{b8}, we obtain
\begin{equation}\label{v8}
\sum_{n\geq0} Q_4^2(n)q^n \equiv \dfrac{{f_2}^3}{f_1} \pmod 4.
\end{equation} Using \eqref{e1} in \eqref{v8}, we obtain
\begin{equation}\label{j8}
\sum_{n\geq0} Q_4^2(n)q^n \equiv \psi(q)f_2  \pmod 4,
\end{equation}
which is the case of $\beta =0$ of \eqref{g7}. Assume that \eqref{g7} holds for some $\beta \geq 0.$ 

Employing \eqref{u10} and \eqref{p8} into \eqref{g7}, we have
\begin{multline}\label{h8}
\sum_{n\geq0} Q_4^2 \left(p^{2\beta} n +\dfrac{5}{12}(p^{2\beta}-1)\right)q^n \\\equiv \Big[ \sum_{k=0}^{{(p-3)/2}}{\it q}^{{(k^2+k)/2}}f\left(q^{{(p^2+(2k+1)p)/2}},q^{{(p^2+(2k+1)p)/2}}\right)+ q^{{(p^2-1)/8}} \psi(q^{p^2}) \Big]\\ \times \Big[ \sum_{\substack{ m={-(p-1)/2} \\ m \ne {(\pm p-1)/6}}}^{m={(p-1)/2}}(-1)^{m}{\it q}^{(3m^2+m)} f\left( {-\it q}^{{2(3p^2+(6m+1)p)/2}},{-\it q}^{2(3p^2-(6m+1)p)/2}\right)\\   
 +(-1)^{{(\pm p-1)/6}}{\it q}^{2(p^2-1)/24}f_{p^2}\Big]  \pmod 4. \end{multline}
Consider the congruence 
 $$3m^2+m + \dfrac{k^2+k}{2} \equiv \dfrac{5(p^2-1)}{24} \pmod p,  $$ which is equivalent to
 \begin{equation}\label{z5}
  (12k+2)^2+6(2m+1)^2\equiv 0 \pmod p.
  \end{equation}
 Since $\left( \dfrac{-6}{p}\right)=-1 $, the only solution of the above congruence is $m=\dfrac{p-1}{2}$ and $k=\dfrac{(\pm p-1)}{6}.$ Therefore, extracting the terms containing $q^{pn+{5(p^2-1)/12}}$ from both sides of \eqref{h8}, dividing by $q^{{5(p^2-1)/12}}$ and replacing $q^p$ by $q,$ we obtain
\begin{equation}\label{g9}
\sum_{n\geq0} Q_4^2 \left(p^{2\beta +1} n +\dfrac{5}{12}(p^{2\beta +2}-1)\right)q^n \equiv \psi(q^p)f_{2p} \pmod 4, \end{equation}
which yields 
\begin{equation}\label{g56}
\sum_{n\geq0} Q_4^2 \left(p^{2\beta+2} n +\dfrac{5}{12}(p^{2\beta+2}-1)\right)q^n \equiv \psi(q)f_2 \pmod 4, \end{equation} 
which is the case $\beta +1$ of \eqref{g7}. Thus, by principle of mathematical induction we complete the proof of \eqref{g7}. 

  Extracting the  coefficients of $q^{pn+j},$ for $1\leq j\leq p-1,$ from both sides of \eqref{g9}, we arrive at  \eqref{s6}.
\end{proof}

\begin{theorem}\label{25} We have
\begin{equation}\label{26}
Q_4^3(2n+1)=0
\end{equation} and 
\begin{equation}\label{27}
Q_4^3(2n)=p_o(n).
\end{equation}
\end{theorem}
\begin{proof}Setting $t=4$ and $s=3$ in \eqref{z2}, we obtain
\begin{equation}\label{28}
\sum_{n\geq0} Q_{4}^{3}(n)q^n=\dfrac{f_2 f_{4}}{f_1 f(q^3,q)}. \end{equation} Using \eqref{e1} in \eqref{28}, we obtain
\begin{equation}\label{29}
\sum_{n\geq0} Q_{4}^{3}(n)q^n=\dfrac{ f_{4}}{f_2}. \end{equation}
Extracting the term involving  $q^{2n+1}$ from \eqref{29},  we arrive at  \eqref{26}.

Again, extracting the term involving  $q^{2n}$  from \eqref{29} and replacing $q^2$ by $q$, we obtain
\begin{equation}\label{30}
\sum_{n\geq0} Q_{4}^{3}(2n)q^n=\dfrac{ f_{2}}{f_1}=\dfrac{(q^2;q^2)_{\infty}}{(q;q)_{\infty}}. \end{equation}
Simplifying \eqref{30} by using elementary $q$-operations, we obtain
\begin{equation}\label{31}
\sum_{n\geq0} Q_{4}^{3}(2n)q^n=\dfrac{1}{(q;q^2)_{\infty}}=\sum_{n\geq0} p_o(n)q^n. \end{equation} Equating the coefficients of $q^n$ on both sides of \eqref{24}, we arrive at \eqref{27}.
\end{proof}

\begin{remark}
Setting $n=2$ in result \eqref{26}, it is easily seen that $Q_4^3(5)=0.$ Again, $Q_4^3(6)=2$ with relevant partition $6,$ and $ 4+2$, and $p_o(3)=2$ with relevant partition $3,$ and $ 1+1+1.$
\end{remark}

\begin{theorem}We have
\begin{equation}\label{32}
 Q_6^2(n)\equiv b_6(n) \pmod2,
\end{equation}
where $b_6(n)$ is the number of partitions of $n$ with no parts divisible by 6.
\begin{equation}\label{33}
 Q_{12}^2(3n+2)\equiv 0 \pmod2.
\end{equation}
\end{theorem}
\begin{proof}
Setting $t=6$ and $s=2$ in  \eqref{z2}, we obtain 
\begin{equation}\label{34}
\sum_{n\geq0} Q_6^2(n)q^n=\dfrac{f_2 f_6}{f_1 f(q^2,q^4)}.
\end{equation} Using \eqref{t8} in \eqref{34}, we obtain
\begin{equation}\label{35}
\sum_{n\geq0} Q_6^2(n)q^n=\dfrac{f_2^2 f_6 f_{12}}{f_1 f_6^2 f_4}.
\end{equation}
 Employing \eqref{t7} in \eqref{35}, we obtain
\begin{equation}\label{36}
\sum_{n\geq0} Q_6^2(n)q^n\equiv \dfrac{ f_6 }{f_1}=\sum_{n\geq0} b_6(n)q^n \pmod2.
\end{equation}
Equating the coefficient of $q^n$ on both sides of \eqref{36}, we arrive at \eqref{32}.

Again,  setting $t=12$ and $s=2$ in \eqref{z2}, we obtain
\begin{equation}\label{37}
\sum_{n\geq0} Q_{12}^2(n)q^n=\dfrac{f_2 f_{12}}{f_1 f(q^2,q^{10})}.
\end{equation}
 Using \eqref{t9} in \eqref{37}, and then simplifying using \eqref{t7}, we have
\begin{equation}\label{39}
\sum_{n\geq0} Q_{12}^2(n)q^n\equiv \dfrac{f_1^3}{f_6}\pmod 2.
\end{equation}
Employing \eqref{y9} in \eqref{39}, we have 
\begin{equation}\label{40}
\sum_{n\geq0} Q_{12}^2(n)q^n=\dfrac{a(q^3)-3qf_9^3}{f_6}\pmod 2.
\end{equation}
Since the right hand side of \eqref{40} has no term involving $q^{3n+2}$, extracting the term $q^{3n+2}$ from \eqref{40}, we arrive at the desired result \eqref{33}.
\end{proof}
\begin{remark}
Setting $n=4$ in \eqref{32}, one can see that $Q_6^2(4)=1 $with relevant partition $3+1$  and $b_6(4)=5$ with relevant partitions $4,\quad 3+1,\quad 2+2,\quad 2+1+1,\quad$ and $ 1+1+1+1$. Again, setting $n=1$ in result \eqref{33}, we observe that $Q_{12}^2(5)=2$ with relevant partitions 5 and $4+1$.
\end{remark}

\section*{Compliance with Ethical Standards}

\textbf{Conflict of interest:} The authors declare that there is no conflict of interest regarding the publication of this article.

\textbf{Human and animal rights:} The author declares that there is no research involving human participants and/
or animals in the contained of this paper.

\end{document}